\documentclass[12pt]{amsart}

\usepackage{amsmath,amsfonts,amssymb,amscd}
\title{An obstruction to asymptotic semistability \\
and approximate critical metrics}
\author{Toshiki Mabuchi}
\address{Department of Mathematics, Graduate School of Science, Osaka University, 
Toyonaka, Osaka, 560-0043 Japan}
\date{}
\begin{document}
\maketitle

\footnotetext{\, To appear in Osaka Journal of Mathematics {\bf 41}(2004).
}

\section{Introduction}

For  a polarized algebraic manifold $(M, L)$ with a K\"ahler metric of constant scalar curvature
in the class $c_1(L)_{\Bbb R}$,
we consider the Kodaira embedding 
$$
\Phi_{|L^m|} : M\; \hookrightarrow \;\Bbb P (V_m), \qquad m \gg 1,
$$
where $V_m := H^0(M, \mathcal{O}(L^m))^*$.
Even when  a linear algebraic group of positive dimension 
acts nontrivially and holomorphically on $M$,
 we shall show that 
the vanishing of an obstruction to asymptotic Chow-semistability allows us to
generalize Donaldson's   construction \cite{Don2} of approximate 
solutions for  equations of critical metrics${}^\dagger$ of Zhang \cite{Z2}.
\footnotetext{${}^\dagger$\, In (2.6) below, $\omega = c_1(L;h)$ is called a {\it critical metric\/}
if $K(q, h)$ is a constant function on $M$.
The same concept was later re-discovered by Luo \cite{Luo} 
(see \cite{M4}).}
This generalization plays a crucial role in our forthcoming paper \cite{M4},  
in which  the asymptotic Chow-stability for $(M, L)$ above will be shown
under the vanishing of the obstruction,
even when $M$ admits a group action as above.

\smallskip
\section{Statement of results}

\medskip
Throughout this paper, we assume that $L$ is  an
 ample holomorphic line bundle over
a connected
projective algebraic manifold
$M$.
Let $n$ and $d$ be respectively the dimension of $M$ and the
degree of  the image $M_m := \Phi_ {|L^m|}
(M )$ in the projective space $\Bbb P (V_m)$ with
$m \gg 1$. Then to this image $M_m$,
we can associate a nonzero element
$\hat{M}_m$ of $W_m := \{\operatorname{Sym}^d (V_m )\}^{\otimes n+1}$ such that 
its natural image $[\hat{M}_m ]$ in 
$\Bbb P (W_m)$ is the Chow point associated to the 
irreducible reduced algebraic cycle
$M_m$ on $\Bbb P(V_m)$.
For the natural action of
$H_m:= \operatorname{SL}(V_m)$ on 
$W_m$
and  also on $\Bbb P ( W_m)$,
the subvariety  $M_m$ of $\Bbb P (V_m)$
 is said to be {\it Chow-stable\/} or 
{\it Chow-semistable\/}, according as  the orbit $H_m\cdot
\hat{M}$ is closed in $W_m$ or 
the
origin of $W_m$ is not in
the closure of $H_m\cdot \hat{M}$ in $W_m$.
Fix an increasing sequence 
$$
  m(1) < m(2) < m(3) < \cdots < m(k) < \cdots 
\leqno{(2.1)}
$$ 
of positive integers $m(k)$. For this sequence,  we say that $(M, L)$ is
{\it asymptotically Chow-stable} or {\it asymptotically Chow-semistable}, 
according as for some
$k_0 \gg 1$,  the subvariety $M_{m(k)}$ 
of $\Bbb P(V_{m(k)})$ is  Chow-stable or Chow-semistable for all 
$k \geq k_0$. 

\smallskip
Let $\operatorname{Aut}^0(M)$ denote the
identity component of the group of all holomorphic automorphisms 
of $M$.    Then the maximal connected linear algebraic subgroup  $G$ of
$\operatorname{Aut}^0(M)$ is the identity component of the kernel
 of the Jacobi homomorphism 
$$ 
\alpha_M : \operatorname{Aut}^0(M) 
\to \operatorname{Aut}^0(\operatorname{Alb}(M)),
\quad \text{(cf. \cite{F1}).}
$$ 
For the maximal algebraic torus $Z$ in the center of $G$, 
we consider the Lie subalgebra $\frak z$ of $H^0(M, \mathcal{O}(T^{1,0}M))$
associated to the  Lie subgroup $Z$ of $\operatorname{Aut}^0(M)$.
For the isotropy subgroup, denoted by $\tilde{S}_m$, 
of $H_m$ at the point $[\hat{M}_m] \in \Bbb P
(W_m)$, we have a natural isogeny 
$$
\iota_m : \tilde{S}_m \to S_m,
$$
where $S_m$ is an algebraic subgroup of $G$. 
For $Z_m := \iota_m^{-1}(Z)$, we have a $Z_m$-action on $M$
naturally induced by the $Z$-action on $M$.
Since the $Z$-action on $M$ 
is liftable to a holomorphic bundle action on $L$ (see for instance \cite{FM1}),  
the restriction of $\iota_m$  to $Z_m$ defines an isogeny of $Z_m$ onto $Z$.
The vector space $V_m$ is viewed as
the line bundle $\mathcal{O}_{\Bbb P(V_m)}(-1)$
with the zero section blown-down to a point,
while the line bundle
$\mathcal{O}_{\Bbb P(V_m)}(-1)$ coincides with $L^{-m}_{}$
when restricted to $M$.
Hence,  the natural $\tilde{S}_m$-action on $V_m$ induces 
a bundle action of $Z_m$ on $L^m$ 
which covers the $Z_m$-action on $M$.
Infinitesimally, each   $X\in\frak z$  
induces a holomorphic vector field
  $X'\in  H^0(L^m, \mathcal{O}(T^{1,0}L^m))$ on $L^m$.
Since  the $\Bbb C^*$-bundle $ L \setminus \{0\}$ associated to $L$
is an $m$-fold unramified covering of the $\Bbb C^*$-bundle $L^m\setminus \{0\}$,
the restriction of  $X'$ to 
$L^m\setminus \{0\}$ naturally induces  a holomorphic vector field 
$X''$ on
$L\setminus \{0\}$. Since $X''$ extends to a holomorphic vector field on $L$,
the mapping $X \mapsto X''$ defines  inclusions
$$
\rho_m : \; \frak z\;  \hookrightarrow\; H^0(L, \mathcal{O}(T^{1,0}L)),
\qquad \text{$m =1,2,\dots$,}
\leqno{(2.2)}
$$
inducing lifts, from $M$ to $L$, of vector fields in $\frak z$. 
For a sequence as in (2.1),  we say that {\it the isotropy actions for $(M, L)$
are stable\/} if there exists an integer $k_0\gg 1$ such that 
$$
\rho_{m(k)} = \rho_{m(k_0)}, \qquad \text{for all  $k \geq k_0$.}
\leqno{(2.3)}
$$
For the maximal compact subgroup $(Z_m)_c$  of $Z_m$,
take a $(Z_m)_c$-invariant Hermitian metric $\lambda$ for $L^m$.
By the theory of equivariant cohomology (\cite{BV}, \cite{FM2}), 
we define (see \cite{MN}, \cite{M1}):
$$
\; \mathcal{C}\{c_1^{n+1}; L^m \}  (X)\; :=\;
\frac{\sqrt{-1}}{2\pi}\,(n+1)\int_M
\lambda^{-1}(X\lambda)\, c_1(L^m; \lambda )^n,
\quad X\in \frak z,
 \leqno{(2.4)}
$$
where  $X\lambda$ is  as in \cite{M1}, (1.4.1). Then the $\Bbb C$-linear map
$\,\mathcal{C}\{c_1^{n+1}; L^m \} : \frak z \to\Bbb C\,$ which sends each $X \in \frak z$ to
$\,\mathcal{C}\{c_1^{n+1}; L^m \}  (X)\in \Bbb C\,$ is independent of the choice of $h$.
The following gives  an obstruction to asymptotic Chow-semistability 
(see \cite{F2}, \cite{MN}, \cite{MW} for related results):

\medskip
{\bf Theorem A.} \; {\em   For a sequence as in $(2.1)$, assume that $(M,L)$ is
asymptotically Chow-semistable. Then for some $k_0\gg 1$, 
the equality $\mathcal{C}\{ c_1^{n+1}; L^{m(k)}\} =
0$  holds for all $k \geq k_0$.
In particular,  for this sequence,
the isotropy actions for $(M,
L)$  are stable.}

\medskip
The following modification of a result in \cite{FM1}  shows that, as an
obstruction, 
 the stability condition (2.3) is essential, since the
vanishing of (2.4) is straightforward from (2.3).

\medskip
{\bf Theorem B.} \; {\em For  sufficiently large $(n+2)$ distinct integers
$m_k$, $k = 0, 1, \dots,  n+1$, suppose that
$\rho_{m_0} = \rho_{m_1} =  \dots  = \rho_{m_{n+1}}$.
Then $\mathcal{C}\{ c_1^{n+1}; L^{m_k}\} = 0$  for all 
$k$.}

\medskip
If $\dim Z = 0$, by setting $m(k) = k$  in (2.1) for all $k >0$, 
we see that $\rho_m$ are trivial for
all
$m \gg 1$, and consequently (2.3) holds.
Note also that
  Donaldson's  result \cite{Don2}  treating the case
 $\dim G = 0$  depends on his construction
of approximate solutions for
equations of  critical metrics of Zhang \cite{Z2}.
In Theorem C down below, assuming (2.3),
we generalize Donaldson's construction to the case $\dim G >0$.

\medskip
 Put $N_m := \dim_{\Bbb C} V_m -1$.
Let $h$ be a Hermitian  metric for $L$ such that $\omega = c_1(L; h)$ is a 
K\"ahler metric on $M$.  
By the inner product
$$
(\sigma,  \sigma')^{}_{h}\; :=\; \int_M <\sigma, \sigma'>^{}_h \, \omega^n, 
\qquad \sigma, \sigma'\in V^*_m,
\leqno{(2.5)}
$$
on $V_m^* = H^0(M, \mathcal{O}(L^{ m}))$, we choose a unitary basis 
$\{ \sigma^{(m)}_0,  \sigma^{(m)}_1, \dots,
\sigma^{(m)}_{N_m}\}$  
for $V_m^*$.
Here, $<\sigma, \sigma'>_h$ denotes the function on $M$ 
obtained as the the pointwise inner product of the sections $\sigma$, $\sigma'$ by
the Hermitian metric $h^m$ on $L^m$. 
Put
$$
K(q, h):= \frac{n!}{m^n}\sum_{i=0}^{N_m} \| \sigma^{(m)}_i \|^{\,2}_h,
\leqno{(2.6)}
$$
where $\|\sigma\|_h^{\,2} :=<\sigma, \sigma>_h$ for all $\sigma \in V_m^*$,
and we set $ q := 1/m$. We then have the asymptotic expansion of Tian-Zelditch (cf. \cite{T1}, \cite{Z1}) 
for $m \gg 1$: 
$$
K(q, h) \; = \; 1 + a_1(\omega ) q + 
a_2 (\omega ) q^2 + a_3(\omega ) q^3 +  \dots\,,
\leqno{(2.7)}
$$
where $a_i (\omega ) $, $i= 1, 2$, \dots, are smooth functions on $M$. 
Then $a_1(\omega ) =  \sigma_{\omega}/2$ (cf. \cite{Lu})
for the scalar curvature 
$ \sigma_{\omega}$
of $\omega$. Put $C_q := \{m^n c_1 (L)^n [M]/n!\}^{-1} (N_m +1)$.
Then

\medskip
{\bf Theorem C.} \; {\em For a K\"ahler metric $\omega_0$
in the class $c_1(L)_{\Bbb R}$ of constant scalar curvature,
choose a Hermitian metric $h_0$ for $L$ such that $\omega_0 = c_1(L, h_0)$.
For a sequence as in $(2.1)$,  assume that
 the isotropy actions for $(M, L)$ are stable,
i.e., $(2.3)$ holds.  Put $q = 1/m(k)$. Then there exists a sequence of 
real-valued smooth functions $\varphi_k$, $k=1,2,\dots,$ 
on $M$ such that 
$h(\ell ):=  h_0 \exp (-\Sigma_{k=1}^{\ell} q^k \varphi_k )$
satisfies 
$K(q, h(\ell )) - C_q = O(q^{\ell +2})$ for each nonnegative integer $\ell$.
}

\medskip\noindent
The last equality $K(q, h(\ell )) - C_q = O(q^{\ell +2})$  means that
there exist a positive real constant $A = A_{\ell}$ independent
of $q$ such that   $\|\, K(q, h(\ell )) - C_q\,\|_{C^0(M)} \leq A_{\ell}\, q^{\ell + 2}$ 
 for all $0 \leq  q \leq 1$ on $M$. By \cite{Z1}, 
for every nonnegative integer $j$,
 a choice of a larger constant $A= A_{j,\ell}>0 $ 
keeps  Theorem C  still valid 
even if  $C^0(M)$-norm is replaced by
$C^j(M)$-norm.

\section{An obstruction to asymptotic semistability}

The purpose of this section is to prove Theorems A and B. Fix a sequence as in (2.1),
and in this section, 
any kind of stability is considered with respect to this
sequence.

\medskip
{\rm Proof of Theorem A}:
Assume that $(M, L)$ is asymptotically Chow-semistable, i.e., for some
$k_0 \gg 1$,
the subvariety $M_{m(k)}$ 
of $\Bbb P (V_{m(k)})$ is Chow-semistable for all $k \geq k_0$.
Then the isotropy representation of $Z_{m(k)}$ on the line 
$\Bbb C \cdot \hat{M}_{m(k)}$ is trivial (cf. \cite{F2}, \cite{MN}) for $k \geq k_0$, and hence by 
  \cite{MN}, (3.5) (cf. \cite{MW}; \cite{Z2}, (1.5)),
we obtain the required equality
$$
\mathcal{C} \{c_1^{n+1}; L^{m(k)}\}(X)\;  = \; 0,  \qquad X \in \frak z,
\leqno{(3.1)}
$$
for all $k \geq k_0$.   For $\lambda$ in (2.4), by setting $h:= \lambda^{1/m}$, we have 
a Hermitian metric $h$ for $L$.
Put $\chi_m := \mathcal{C}\{c_1^{n+1}, L^{m}\}/m^{n+1}\, $ for positive integers $m$. 
Then by the Leibniz rule,
$$
\chi_m (X) \; =\; \frac{\sqrt{-1}}{2\pi}\, (n+1)
\int_M h^{-1}(Xh)_{\rho_m}\, c_1(L; h)^n, \qquad X \in \frak z,
\leqno{(3.2)}
$$
where the complexified action $(Xh)_{\rho_m}$ of $X$ on $h$   as in \cite{M1}, (1.4.1),
 is taken via the lifting $\rho_m$ in (2.2). Then by (3.1),
$$
\chi_{m(k_0)}\; =\;\chi_{m(k_0+1)}\; =\; \dots \dots\; =\; \chi_{m(k)}
\; =\; \dots\dots,
$$
and since lifts in (2.2), from $M$ to $L$,  of holomorphic
vector fields in $\frak z$  are completely characterized by
$\chi_m$ (cf. \cite{FM1}), we obtain (2.3), as required.
\qed

\medskip
{\rm Proof of Theorem B}: \; For $q :=  \text{l.c.m}\{m_k; k= 0,1,\dots, n+1\}$,
we take a $q$-fold unramified cover $\nu : \tilde{Z} \to Z$ between algebraic tori. Then
the $Z$-action on $M$ naturally induces a $\tilde{Z}$-action on $M$ via this covering.  
Since $\nu$ factors 
through $Z_{m_k}$,
the lift,  from $M$ to $L^{m_k}$, of the $Z_{m_k}$-action
naturally induces a lift, from $M$ to $L^{m_k}$,  of the $\tilde{Z}$-action.
The assumption
 $$
\rho_{m_0} = \rho_{m_1} = \dots = \rho_{m_{n+1}} 
\leqno{(3.3)}
$$ 
shows that
the lifts, from $M$ to $L^{m_k}$,  $k$ = 0,1,\dots, $n+1$,  of the $\tilde{Z}$-action
come from the same infinitesimal action of $\frak z$ as vector fields on $L$.  For brevity, the
common
$\rho_{m_k}$ in (3.3)  will be denoted just by $\rho$. Then the proof of \cite{FM0}, Theorem 5.1, is
valid also in our case, and the formula in the theorem holds.
By $Z_{m_k} \subset \operatorname{SL}(V_{m_k})$ and by its contragredient 
representation, 
the $\tilde{Z}$-action on $V_{m_k}^* = H^0(M, \mathcal{O}(L^{m_k}))$ 
comes from an algebraic group homomorphism:
$\tilde{Z} \to \operatorname{SL}(V_{m_k}^*)$. Hence,
by  the notation in (3.2) above,
$\int_M h^{-1}(Xh)_{\rho} \, c_1(L; h)^n  = 0$ for all $X \in \frak z$,
i.e., $\mathcal{C}\{c_1^{n+1}; L^{m_k}\}  = 0$ for all $k$, as required.
\qed

\section{Proof of Theorem C}

\smallskip
Throughout this section, we assume that the first Chern class $c_1(L)_{\Bbb R}$ admits
a K\"ahler metric of constant scalar curvature. Then  a result of Lichn\'erowicz \cite{Lich}
(see also \cite{kob1}) shows that $G$ is a reductive algebraic group, 
and consequently the identity component of the center of $G$ coincides with $Z$ 
in the introduction. Let $K$ be a maximal compact subgroup of $G$.
Then  the maximal compact subgroup $Z_c$ of $Z$ satisfies
$$
Z_c \; \subset \; K.
\leqno{(4.1)}
$$
For an arbitrary  $K$-invariant K\"ahler metric $\omega$
 on $M$ in the class $c_1(L)_{\Bbb R}$, we  write
$\omega$  as the Chern form $c_1(L;h)$ 
for some Hermitian metric $h$ for $L$. 
Let $\Psi (q, \omega)$ denote 
the power series in $q$ given by the right-hand side of (2.7). 
Then
$$
\int_M  \{ \Psi (q, \omega) - C_q \} \,\omega^n 
\; =\; \int_M\left \{ - C_q + \frac{n!}{m^n}\sum_{i=0}^{N_m} \| \sigma^{(m)}_i \|^{\,2}_h \right \}\omega^n
\; =\; 0.
\leqno{(4.2)}
$$
Let $h_0$ be a Hermitian metric for $L$ such that $\omega_0 := c_1(L; h_0)$ is a K\"ahler metric
of constant scalar curvature on $M$. We write
$$
\omega_0 \; =\; \frac{\sqrt{-1}}{2\pi}\;\sum_{\alpha, \beta}\; g_{\alpha\bar{\beta}} 
dz^{\alpha}\wedge d
z^{\bar{\beta}},
$$
for a system $(z^1, z^2, \dots , z^n)$ of holomorphic local coordinates on $M$.
In view of \cite{Lich} (see also \cite{kob1}), replacing $\omega_0$ by $g^*\omega_0$ for some 
$g \in G$ if necessary, we may assume that $\omega_0$ is $K$-invariant. 
Let $D_0$ be the Lichn\'erowicz operator, as defined in \cite{Ca}, (2.1), for the K\"ahler manifold
$(M,
\omega_0 )$. Since $\omega_0$ has a constant scalar curvature, $D_0$ is a real operator.
Let $\mathcal{F}$ denote the space of all real-valued smooth $K$-invariant functions $\varphi$ 
such that $\int_M \varphi \omega_0^{\,n} = 0$.
Since the operator $D_0$ preserves the space $\mathcal{F}$, we write $D_0$
as an operator $D_0 :
\mathcal{F} \to \mathcal{F}$, and the kernel in $\mathcal{F}$ of this operator
will be denoted by $\operatorname{Ker} D_0$.  
Let $\frak z_c$ denote the Lie subalgebra of $\frak z$
corresponding to the maximal compact subgroup $Z_c$ of $Z$.
Then 
$$
\gamma : \operatorname{Ker} D_0 \; \cong \; \frak z_c,
\qquad \eta \leftrightarrow 
\gamma (\eta ) := \operatorname{grad}_{\omega_0}^{\Bbb C}\eta,
\leqno{(4.3)}
$$
where $\operatorname{grad}_{\omega_0}^{\Bbb C}\eta := 
 (1/\sqrt{-1}) \Sigma\, g^{\bar{\beta}\alpha}\eta_{\bar{\beta}}
\partial /\partial z^{\alpha}$ denotes the complex gradient
of $\eta$ with respect to $\omega_0$.
We then consider the orthogonal
projection
$$
P\;  : \; \mathcal{F} \; (=\operatorname{Ker} D_0 \oplus \operatorname{Ker} D_0^{\perp}) \;\to\;
\operatorname{Ker} D_0.
$$ 
Starting from $h (0) = h_0$ and $\omega (0):= \omega_0$, we inductively
define a Hermitian metric 
$h(k )$ for $L$, and a K\"ahler metric $\omega (k ) := c_1(L; h(k))$,
called the {\it  $k$-approximate solution\/}, by
\begin{align*}
h( k) \; &= \; h (k-1) \, \exp (- q^k \varphi_k ), \qquad\qquad k = 1,2,\dots, \\
\omega (k) \; &=\; \omega (k-1)\, +\, \frac{\sqrt{-1}}{2\pi}
\,q^k\,\partial\bar{\partial}\varphi_k , \qquad k = 1,2,\dots, 
\end{align*}
for a suitable function $\varphi_k \in \operatorname{Ker} D_0^{\perp}$, 
where   we require $h (k)$ 
to satisfy 
$K (q, h (k)) - C_q = O(q^{k+2})$. In other words,  by (4.2), each $\omega (k)$ 
is required to satisfy the following conditions:
\begin{align*}
 (1-P)\{\Psi (q, \omega (k)) - C_q \} \; &\equiv  \; 0,   \quad \text{modulo $q^{k+2}$},  \tag{4.4} \\
 P \{\Psi (q, \omega (k) ) - C_q \} \; &\equiv \; 0,   \quad \text{modulo $q^{k+2}$}. \tag{4.5}
\end{align*}

If $k =0$, then $\omega (0) = \omega_0$, and by \cite{Lu}, both (4.4) and (4.5) hold for 
$k = 0$. Hence, let $\ell \geq 1$ and assume  (4.4) and (4.5)  for $ k = \ell -1$.
It then suffices to find $\varphi_{\ell} \in \operatorname{Ker} D_0^{\perp}$
satisfying both (4.4) and (4.5) for $k = \ell$.
Put
$$
\Phi (q, \varphi )\;  := \; (1-P)\left \{\,\Psi \left (q, \,\omega (\ell -1) + (\sqrt{-1}/2\pi )\,q^{\ell} 
\,\partial\bar{\partial}\varphi 
\right ) -  C_q\right \}, \qquad \varphi \in \operatorname{Ker} D_0^{\perp}.
$$
Then by (4.4) applied to $k = \ell - 1$, we have 
$\Phi (q, 0) \equiv\, u_{\ell}\, q^{\ell + 1}$ modulo $q^{\ell + 2}$,
where $u_{\ell}$ is a  function in $\operatorname{Ker} D_0^{\perp}$.
Since $2 \pi \omega (\ell -1) = 2 \pi \omega_0 + \sqrt{-1}\;\Sigma_{k=1}^{\ell - 1}\, q^k
\partial\bar{\partial}\varphi_k$, we have $\omega (\ell -1) = \omega_0$ at $q  = 0$.
Since the scalar curvature of $\omega_0$ is constant,
 the variation formula for the scalar curvature (see for instance
\cite{Ca}, (2.5);  \cite{Don2}) shows that
$$
 \Phi (q, \varphi_{\ell} ) \; \equiv \; \Phi (q, 0 ) 
- q^{\ell +1} (D_0 \varphi_{\ell}/2) \; \equiv \; ( 2u_{\ell} - D_0 \varphi_{\ell})\,(q^{\ell +1}/2), 
$$
modulo $q^{\ell +2}$. Since $u_{\ell}$ is in $\operatorname{Ker} D_0^{\perp}$,
there exists a unique $\varphi_{\ell}\in \operatorname{Ker} D_0^{\perp}$ 
such that $ 2u_{\ell} = D_0 \varphi_{\ell}$ on $M$. Fixing such $\varphi_{\ell}$,
we obtain $h (\ell )$ and $\omega (\ell )$.
Thus (4.4) is true for $k = \ell$. 

\smallskip
Now, we have only to show that (4.5) is true for $k = \ell$.  
Before checking this, we give some preliminary remarks.
Note that $C_q = 1 + O(q)$.
Moreover, by (2.7), $\Psi (q, \omega ) = 1 + q\{ a_1 (\omega ) + a_2 (\omega ) q + \dots \}$,
and hence
\begin{align*}
&\Psi (q, \omega (\ell ) )-C_q \; =\; \Psi \left (\,q,\,\omega (\ell -1) + (\sqrt{-1}/2\pi )\, q^{\ell}
\partial\bar{\partial}\varphi_{\ell} \right )  - C_q \\
&\equiv \; \Psi (q, \omega (\ell -1)  ) - C_q \; \equiv\; 0,
\qquad \quad \text{modulo $q^{\ell + 1}$}.
\end{align*}
By \cite{MFK}, p. 35,  the $G$-action on $M$ is liftable to a bundle action of $G$
on the real line bundle $(L\cdot \bar{L})^{1/2} = (L^{m}\cdot\bar{L}^{m})^{1/2m}$.
Then the induced $K$-action on $(L\cdot \bar{L})^{1/2}$ is unique,
because liftings, from $M$ to $L^m$, of the $G$-action differ only by scalar multiplications of $L^m$
by  characters of $Z$.
In this sense, $h (\ell )$ is $K$-invariant. 
Put $r:= \dim_{\Bbb C} Z$. 
Then we can write 
$Z_m \; = \; \Bbb G_m^r \; =\; \{\, t = (t_1, t_2, \dots , t_r) \in (\Bbb C^*)^r\,\}$.
By the natural inclusion
$$
\psi_m \;:\; Z_m\;  \hookrightarrow\;  H_m = \operatorname{SL}(V_m),
$$
we can choose a unitary basis $\{\tau_0, \tau_1, \dots, \tau_{N_{m}}\}$ for 
$(V^*_{m}, (\; , \, )^{}_{h (\ell ) })$ (cf. (2.5))
such that, for some integers $\alpha_{ij}$ with $\Sigma_i\, \alpha_{ij} = 0$,
the contragredient representation $\psi_m^*$ of $\psi_m$ is given by
$$
\psi^*_m (t) \,\tau_i \; =\; \left (\prod_{j=1}^r \, t_j^{\alpha_{ij}} \right ) \tau_i,
\qquad i = 0,1, \dots, N_m, 
$$
for all $t \in (\Bbb C^*)^r = Z_m$.
Now by (2.3),
for some $\rho: \frak z \hookrightarrow H^0(L, \mathcal{O}(T^{1,0}L))$,
we can write $ \rho_{m(k)} = \rho$ for all $k \geq k_0$. 
Consider the K\"ahler metric 
$\omega_m := c_1(L; h_m)$ on $M$ in the clasas $c_1(L)_{\Bbb R}$,
where $h_m := ( |\tau_0|^2 + |\tau_1 |^2 + \dots + |\tau_{N_m}|^2 )^{-1/m}$.
From now on, let $m = m(k)$, where $k$ is running through all integers $ \geq k_0$.
Put $X_j  := t_j \partial/\partial t_j$. 
Then  $\{X_1, X_2, \dots, X_r\}$ forms a $\Bbb C$-basis for the Lie algebra $\frak z$  
 such that, using the notation as in (3.2), we have
$$
h_m^{-1}(X_j h_m)_{\rho}\; =\; -\,  \frac{\Sigma_i\,  \alpha_{ij} |\tau_i |^2 }{m\,\Sigma_i\,   |\tau_i |^2 },
\quad  1 \leq j \leq r,\quad \text{for $m = m(k)$ with $k \geq k_0$,}
\leqno{(4.6)}
$$
where in the numerator and the denominator, the sum is taken over all 
integers $i$  such that $0\leq i \leq N_m$. From (2.3) and Theorem B, using the notation as in
(3.2),  we obtain
$$
\int_M   h(\ell  )^{-1}
( X_j h(\ell  ))_{\rho}\, \omega (\ell )^n\; = \; 0, \qquad  1 \leq
j
\leq r.
\leqno{(4.7)}
$$
By $\int_M h_0^{-1}(X_j h_0)_{\rho}\,\omega_0^{\,n}/\int_M \omega_0^{\,n} = 0$, we have
$\eta_j :=  h_0^{-1}(X_j h_0)_{\rho} \in \operatorname{Ker} D_0$.
Then $\gamma ( \eta_j ) = \sqrt{-1}\, X_j$.  Hence
$\{ \eta_1, \eta_2, \dots,
\eta_r \}$ is an $\Bbb R$-basis for $\operatorname{Ker} D_0$.
Since $\Psi (q, \omega (\ell )) \equiv C_q $ modulo $q^{\ell + 1}$, it follows that
$$
 - C_q + \frac{n!}{m^n}\sum_{i=0}^{N_m} \| \tau_i \|^{\,2}_{h (\ell ) } \; \equiv \;  v_{\ell}\,q^{\ell + 1}
\leqno{(4.8)}
$$
modulo $q^{\ell + 2}$ for some $v_{\ell}\in \operatorname{Ker} D_0$,
because (4.4) is true for $k = \ell$. 
In view of (4.2), (4.6),
$h_m - h_0 = O(q)$ and $\omega (\ell ) - \omega_0 = O(q)$,  we see from (4.8) that,
 modulo
$q^{\ell + 2}$,
\begin{align*}
q^{\ell +1}\int_M  \eta_j \,v_{\ell} \, \omega_0^{\,n} \; &\equiv \; 
\int_M \eta_j \left (  - C_q + \frac{n!}{m^n}\sum_{i=0}^{N_m} \| \tau_i \|^2_{h (\ell )}
\right ) \{\omega (\ell  )\}^n\\
\; &\equiv \; \int_M h_0^{-1}(X_j h_0)_{\rho} 
\left (  - C_q + \frac{n!}{m^n}\sum_{i=0}^{N_m} \| \tau_i \|^2_{h (\ell )}
\right ) \{\omega (\ell  )\}^n \\
\; &\equiv \; \int_M h_m^{-1}(X_j h_m)_{\rho} 
\left (  - C_q + \frac{n!}{m^n}\sum_{i=0}^{N_m} \| \tau_i \|^2_{h (\ell )}
\right ) \{\omega (\ell  )\}^n \\
\; & \equiv \; 
\int_M \frac{\Sigma_i\,  \alpha_{ij} \|\tau_i \|^2_{h (\ell ) } }{m \,\Sigma_i\,   \|\tau_i \|^2_{h(\ell )} }
\left (   C_q - \frac{n!}{m^n}\sum_{i=0}^{N_m} \| \tau_i \|^2_{h(\ell )}
\right ) \{\omega (\ell  )\}^n,
\end{align*}
Since
 $\Sigma_{i} \alpha_{ij} = 0$ for all $j$,
we obtain, modulo $q^{\ell + 2}$,
\begin{align*}
q^{\ell +1}\int_M  \eta_j\,v_{\ell} \, \omega_0^{\,n} \;  &\equiv \; 
 C_q \int_M \frac{\Sigma_i\,  \alpha_{ij} \|\tau_i \|^2_{h(\ell )} }{m\, \Sigma_i\,   \|\tau_i \|^2_{h(\ell )} }
\{\omega (\ell  )\}^n\; \equiv \;   C_q \int_M  h_m^{-1}(X_j h_m )_{\rho}\{\omega (\ell  )\}^n \\
& \equiv \; C_q \int_M \left \{  h_m^{-1}(X_j h_m )_{\rho} -
 h(\ell )^{-1} (X_j h(\ell ))_{\rho}\right \} \{\omega (\ell  )\}^n,
\end{align*}
where  the  equivalence just above follows from (4.7). The last integrand is rewritten as
\begin{align*}
 &h_m^{-1}(X_j h_m )_{\rho} -
 h (\ell )^{-1} (X_j h (\ell ))_{\rho} \, = \,
X_j \log \{ h_m/  h (\ell )  \} \, =\, -\,\frac{1}{m}\,X_j \log
\left ( \frac{n!}{m^n}\sum_{i=0}^{N_m} \| \tau_i \|^2_{h(\ell )} \right )\\
&\equiv \; -\,q\,X_j \log (C_q + v_{\ell}\, q^{\ell +1} ) \; \equiv 
\; -\,C_q^{\,-1} (X_j v_{\ell} ) q^{\ell
+2}\; \equiv \; 0,
\qquad \text{mod $q^{\ell + 2}$.}
\end{align*}
Therefore,
$\int_M  \eta_j v_{\ell} \, \omega_0^{\,n}   =  0$ 
for all $j$. From $v_{\ell}\in \operatorname{Ker} D_0$, it now follows that $v_{\ell} = 0$.
This shows that (4.5) is true for $k = \ell$, as required.
\qed

\section{Concludung remarks}

As in Donaldson's work \cite{Don2}, the construction of approximate solutions 
in Threorem C is a crucial step to the approach  of the stability problem
for a polarized algebraic manifold  with a K\"ahler metric of constant 
scalar curvature.
Actually, in a forthcoming paper \cite{M4}, this construction allows us to prove the following:

\medskip\noindent
{\bf Theorem.}
 {\em For a sequence as in $(2.1)$,  assume that
 the isotropy actions for $(M, L)$ are stable.
Assume further that $c_1(L)_{\Bbb R}$ 
admits a K\"ahler metric of constant scalar curvature. Then for this sequence,
$(M, L)$ is asymptotically Chow-stable.}

\medskip
Moreover, if a sequence (2.1) exists in such a way that (2.3) holds,
then the same argument as in the case $\dim G = 0$ (cf. \cite{Don2}) is applied,
and we can also show the uniquness, modulo the action of $G$, of the K\"ahler metrics
 of constant scalar curvature in the polarization class $c_1(L)_{\Bbb R}$.
We finally remark that, if $\dim G = 0$, the asymptotic Chow-stability
 implies the asymptotic stability in the sense of
Hilbert schemes (cf. \cite{MFK}, p.215).
Hence the result of Donaldson \cite{Don2} follows from the theorem just above.

\end{document}